\documentclass{article}

\usepackage[latin1]{inputenc}
\usepackage{t1enc}
\usepackage{amsmath, amssymb, amsfonts, amsthm, amsopn,epsfig}
\usepackage[none]{hyphenat}

\theoremstyle{plain}
\newtheorem{theorem}{Theorem}
\newtheorem{prop}[theorem]{Proposition}

\newtheorem{conjecture}[theorem]{Conjecture}
\newtheorem{OQ}[theorem]{Question}

\theoremstyle{definition}

\title{Decompositions of the Boolean lattice into rank-symmetric chains}
\author{Istv\'{a}n Tomon\\ University of Cambridge}

\begin{document}

\maketitle

\begin{abstract}
The Boolean lattice $2^{[n]}$ is the power set of $[n]$ ordered by inclusion. A chain $c_{0}\subset...\subset c_{k}$ in $2^{[n]}$ is rank-symmetric, if $|c_{i}|+|c_{k-i}|=n$ for $i=0,...,k$; and it is symmetric, if $|c_{i}|=(n-k)/2+i$. We show that there exist a bijection $$p: [n]^{(\geq n/2)}\rightarrow [n]^{(\leq n/2)}$$ and a partial ordering $<$ on $[n]^{(\geq n/2)}$ satisfying the following properties:
\begin{itemize}
  \item $\subset$ is an extension of $<$ on $[n]^{(\geq n/2)}$;
  \item if $C\subset [n]^{(\geq n/2)}$ is a chain with respect to $<$, then $p(C)\cup C$ is a rank-symmetric chain in $2^{[n]}$, where $p(C)=\{p(x): x\in C\}$;
  \item the poset $([n]^{(\geq n/2)},<)$ has the so called normalized matching property.
\end{itemize}

We show two applications of this result.

A conjecture of F\"{u}redi asks if $2^{[n]}$ can be partitioned into $\binom{n}{\lfloor n/2\rfloor}$ chains such that the size of any two chains differ by at most 1. We prove an asymptotic version of this conjecture with the additional condition that every chain in the partition is rank-symmetric: $2^{[n]}$ can be partitioned into $\binom{n}{\lfloor n/2\rfloor}$ rank-symmetric chains, each of size  $\Theta(\sqrt{n})$.

Our second application gives a lower bound for the number of symmetric chain partitions of $2^{[n]}$. We show that $2^{[n]}$ has at least $2^{\Omega(2^{n}\log n/\sqrt{n})}$ symmetric chain partitions.
\end{abstract}

\section{Introduction}

Let us introduce the main definitions and notation used throughout the paper. The notation is mostly standard and can be found in \cite{anderson}, for example.

The \emph{Boolean lattice} $2^{[n]}$ is the power set of $[n]=\{1,...,n\}$ ordered by inclusion. The family of $k$ element subsets of $[n]$ is denoted by $[n]^{(k)}$; also $[n]^{(\leq k)}=\{x\in 2^{[n]}: |x|\leq k\}$, and define $[n]^{\geq k}$ similarly.

A \textit{chain} in a poset is a subset of pairwise comparable elements. A chain  $C\subset 2^{[n]}$ with elements $c_{0}\subset ...\subset c_{k}$ is \emph{rank-symmetric}, if ${|c_{i}|+|c_{k-i}|=n}$, and $C$ is \emph{skipless}, if $|c_{i}|=|c_{0}|+i$ for $i=0,...,k$. A chain is \emph{symmetric}, if it is rank-symmetric and skipless.

A poset $P$ is \emph{graded} if there exists a partition of its elements into subsets $A_{0},A_{1},\ldots,A_{n}$ such that $A_{0}$ is the set of minimal elements, and whenever $x\in A_{i}$ and $x<y$ with no $x<z<y$, then $y\in A_{i+1}$. If there exists such a partition, then it is unique and $A_{0},A_{1},\ldots,A_{n}$ are the levels of $P$. If $x\in A_{i}$, then the \emph{rank} of $x$ is $i$ and it is denoted by $rk(x)$.

A graded poset $P$ is \emph{unimodal} if there exists $0\leq m\leq n$ such that $|A_{0}|\leq\ldots\leq |A_{m}|$ and $|A_{m}|\geq |A_{m+1}|\geq\ldots\geq |A_{n}|$. Also, $P$ is \emph{rank-symmetric}, if we have $|A_{i}|=|A_{n-i}|$ for $i=0,\ldots,n$.

A bipartite graph $G=(A,B,E)$ is a \emph{normalized matching graph}, if for any $X\subset A$ we have
$$\frac{|X|}{|A|}\leq \frac{|\Gamma(X)|}{|B|},$$
where $\Gamma(X)$ is the set of neighbours of $A$ in $B$. A graded poset $P$ is a \emph{normalized matching poset}, or satisfies \emph{the normalized matching property}, if  the bipartite graph induced any two levels $A_{i}$ and $A_{j}$ is a normalized matching graph for $i,j\in \{0,...,n\}$, $i\neq j$.

It is easy to show that $2^{[n]}$ is a rank-symmetric, unimodal normalized matching poset.

\bigskip

By the classical theorem of Sperner \cite{sperner}, the minimum number of chains $2^{[n]}$ can be partitioned into is $\binom{n}{\lfloor n/2\rfloor}$. Also, Brujin et al. \cite{BTK}  showed that $2^{[n]}$ admits a chain decomposition into $\binom{n}{\lfloor n/2\rfloor }$  symmetric chains; later Griggs \cite{symchain} extended their result by proving that every rank-symmetric, unimodal normalized matching poset has a symmetric chain decomposition.

The aim of this paper is to show the existence of a relatively "rich" subposet of $(2^{[n]},\subset)$, where every chain corresponds to a rank-symmetric chain in $2^{[n]}$, and every chain partition corresponds to a chain partition of $2^{[n]}$ into rank-symmetric chains. Here, rich means that the poset has the normalized matching property. Our main result is the following.

\begin{theorem}\label{mainthm2}
For every positive integer $n$ there exist a bijection $$p: [n]^{(\geq n/2)}\rightarrow [n]^{(\leq n/2)}$$ and a partial ordering $<$ on $[n]^{(\geq n/2)}$ such that
\begin{description}
  \item[(i)] $\subset$ is an extension of $<$;
  \item[(ii)] if $C\subset[n]^{(\geq n/2)}$ is a chain with respect to $<$, then $C\cup p(C)$ is a rank-symmetric chain in $2^{[n]}$, where $p(C)=\{p(x): x\in C\}$;
  \item[(iii)] the poset $([n]^{(\geq n/2)},<)$ satisfies the normalized matching property.
\end{description}
\end{theorem}

We also show two applications of this theorem.

The first application is the following  rank-symmetric variant of a problem of F\"{u}redi.
Note that in symmetric chain decomposition of $2^{[n]}$, we must have $\binom{n}{i}-\binom{n}{i-1}$ chains of size $n+1-2i$ for $i=0,...,\lfloor n/2\rfloor $. In other words, this chain decomposition contains small an large chains as well. In 1985, F\"{u}redi \cite{furedi} proposed the following conjecture.

\begin{conjecture}\label{conjfuredi}
Let $n$ be a positive integer. The Boolean lattice $2^{[n]}$ can be partitioned into $\binom{n}{\lfloor n/2\rfloor}$ chains such that the size of any two chains differ by at most 1.
\end{conjecture}

If the conjecture is true, it means that $2^{[n]}$ admits a chain partition into $\binom{n}{\lfloor n/2\rfloor}$ chains such that the size of each chain is $\approx\sqrt{\pi/2}\sqrt{n}$. While this conjecture is still open, there are a few partial results. Hsu, Logan, Shahriari and Towse \cite{hsu} proved that there exists a chain partition into $\binom{n}{\lfloor n/2\rfloor}$ skipless chains such that the size of each chain is at least $\sqrt{n}/2+O(1)$. Also, the author of this paper \cite{me} proved that there is a chain partition of $2^{[n]}$ into $\binom{n}{\lfloor n/2\rfloor}$ chains such that the size of each chain is between $0.8\sqrt{n}$ and $13\sqrt{n}$. We also proposed the following rank-symmetric version of Conjecture \ref{conjfuredi}.

\begin{conjecture}\label{conjsym}
Let $n$ be a positive integer. The Boolean lattice $2^{[n]}$ can be partitioned into $\binom{n}{\lfloor n/2\rfloor}$ rank-symmetric chains such that the size of any two chains differ in at most 2.
\end{conjecture}

 To demonstrate the difficulty of the problem, we challenge the reader to think about the following much weaker result: if $h$ is fixed and $n$ is sufficiently large, then $2^{[n]}$ can be partitioned into rank-symmetric chains, each of size at least $h$. While this problem is not too hard in the case we do not demand our chains to be rank-symmetric, we have to overcome extra obstacles in the rank-symmetric case.

 In this paper, we prove the following result concerning Conjecture \ref{conjsym}.

\begin{theorem}\label{mainthm}
Let
$$\alpha=\sqrt{2}\displaystyle\sum_{k=2}^{\infty}\frac{\sqrt{\log k}-\sqrt{\log (k-1)}}{k}\approx 0.8482.$$
For any $\epsilon>0$ there exists an  $N_{\epsilon}$ such that if $n>N_{\epsilon}$, then $2^{[n]}$ can be partitioned into $\binom{n}{\lfloor n/2\rfloor}$ rank-symmetric chains, all of them of size between $(\alpha-\epsilon)\sqrt{n}$ and $O(\sqrt{n}/\epsilon)$.

In particular, if $n$ is sufficiently large, then there is a chain partition of $2^{[n]}$ into $\binom{n}{\lfloor n/2\rfloor}$ rank-symmetric chains such that the size of each chain is between $0.8\sqrt{n}$ and $13\sqrt{n}$.
\end{theorem}

We note that this is exactly the same result as Theorem 1.5 in \cite{me}, but with the additional condition that our chains are rank-symmetric.

\bigskip

With our second application, we show that $2^{[n]}$ has a lot of symmetric chain decompositions. We prove the following theorem.

\begin{theorem}\label{nbofpartitions}
The Boolean lattice $2^{[n]}$ has at least
$$2^{\Omega(2^{n}\log n/\sqrt{n})}$$
different partitions into $\binom{n}{\lfloor n/2\rfloor}$ symmetric chains.
\end{theorem}

In the proof, we also establish a nontrivial lower bound for the number of matchings in an arbitrary normalized matching graph.

\section{The proof of Theorem \ref{mainthm2}}

In this section, we prove Theorem \ref{mainthm2}. The proof relies on the idea of the classical symmetric chain decomposition of $2^{[n]}$ introduced by De Bruijn, Tengbergen, Kruyswijk \cite{BTK}. We shall briefly define this chain partition.

\bigskip

First, we introduce some notation. If $v$ is an element of a cartesian product with $d$ terms and $i\in [d]$, then $v_{i}$ denotes the $i$th coordinate of $v$.

As usual, let $[i,j]=\{i,i+1,...,j\}$ for $i,j$ integers with $i<j$. Also, for $x\in 2^{[n]}$, let
$$c_{x}(i,j)=|x\cap [i,j]|-|[i,j]\setminus x|.$$ The \emph{signature} of $x$, denoted by $sg(x)$, is an element of $\{0,1,*\}^{n}$ defined as follows.

 If $i\in x$ and there exists an integer $j$ with $i<j\leq n$ such that $c_{x}(i,j)=0$, then $sg(x)_{i}=1$. The smallest such $j$ is the \emph{pair of $i$ in $x$} and is denoted by $pr_{x}(i)$. If there is no such $j$, then $sg(x)_{i}=*$.

 If $i\not\in x$ and there exists an integer $j$ with $1\leq j<i$ such that $c_{x}(j,i)=0$, then $sg(x)_{i}=0$. The largest such $j$ is the \emph{pair of $i$ in $x$} and is denoted by $pr_{x}(i)$. If there is no such $j$, then $sg(x)_{i}=*$.

 Let $*(x)=\{i\in [n]: sg(x)_{i}=*\}$, and define $0(x)$ and $1(x)$ similarly. Also, let $I(x)=\{[i,pr_{x}(i)]: i\in 1(x)\}$ be the set of intervals whose endpoints are the pairs of $x$.

 For example, if $x=\{2,3,4,6,7,10,11\}\in 2^{[12]}$, then $$sg(x)=(*,*,*,\overline{1,0},\overline{1,\overline{1,0},0},*,\overline{1,0}),$$
 so $*(x)=\{1,2,3,10\}$, $1(x)=\{4,6,7,11\}$ and $0(x)=\{5,8,9,12\}$.

 The signature of $x$ can be interpreted in another way as well. Using our earlier example, write $x$ as $)((()(())(()$, where "(" represents the elements in $x$ and ")" represents the elements not in $x$. Then, we have $sg(x)_{i}=*$ if the bracket representing $i$ is invalid, in other words does not have a pair.

 Now we list some of the most important properties of the signature. We shall avoid their proof, as they can be found in \cite{BTK}, for example. Also, the proofs of these properties are very similar to the ones we give in Proposition \ref{properties_of_csg}.

 \begin{prop}\label{properties}
 Let $x\in 2^{[n]}$.
 \begin{description}
   \item[(i)] If $i\in [n]$ and $j=pr_{x}(i)$ exists, then $pr_{x}(j)$ exists as well and $pr_{x}(j)=i$. Also, $sg(x)_{i}=1-sg(x)_{j}$.
   \item[(ii)] $\max(*(x)\setminus x)<\min(x\cap *(x))$.
   \item[(iii)] Any two intervals in $I(x)$ are either disjoint, or one is contained in the other. Also, the intervals are disjoint from $*(x)$.
 \end{description}
 \end{prop}

  \hfill$\Box$

 \bigskip

Define the equivalence relation $\sim$ on $2^{[n]}$ such that $x\sim y$ if $sg(x)=sg(y)$. We show that a complete equivalence class of $\sim$ is a symmetric chain in $2^{[n]}$. Let $x\in 2^{[n]}$ and let $C$ be the $\sim$ equivalence class of $x$. Let $1\leq i_{1}<i_{2}<...<i_{t}\leq n$ be the elements of $*(x)$.  Then one can easily show that the elements of $C$ are the subsets of $[n]$ having the form $$1(x)\cup \{i_{u+1},i_{u+1},...,i_{t}\}$$ for $u=0,...,t$.
 Clearly, $C$ is a chain. Also, using that $|1(x)|=|0(x)|$, we have $|1(x)|=(n-t)/2$. Hence, $$|1(x)\cup \{i_{u+1},...,i_{t}\}|=(n+t)/2-u,$$ and $C$ is truly a symmetric chain. Thus, the equivalence classes of $\sim$ form a symmetric chain decomposition of $2^{[n]}$.

 Now we shall define our bijection $p$. For simplicity, write $Q=[n]^{(\geq n/2)}$. Our bijection $p: Q\rightarrow [n]^{(\leq n/2)}$ is defined with respect to the symmetric chain decomposition described above: for $x\in Q$, let $p(x)$ be the unique element in the $\sim$ equivalence class of $x$ such that $$|x|+|p(x)|=n.$$

Define the relation $<$ on $Q$ as follows. Let $y<x$ if $$p(x)\subset p(y)\subset y\subset x.$$ As $\subset$ is a partial order, $<$ is also a partial order. Hence, $(Q,<)$ is a partially ordered set. Furthermore, $Q$ has the following important property: if $C$ is a chain in $(Q,<)$, then $C\cup p(C)$ is a rank-symmetric chain in $(2^{[n]},\subset)$, where $p(C)=\{p(x): x\in C\}$. Hence, our task is reduced to prove the following theorem.

\begin{theorem}\label{Qnormalized}
The graded poset $(Q,<)$ is a normalized matching poset.
\end{theorem}

The next proposition shall give us a more useable description of $<$.

\begin{prop}\label{comparelemma}
Let $k\geq n/2+1$, $x\in [n]^{(k)}$, and let $i_{1}<...<i_{t}$ be the elements of $x\cap *(x)$. Also, let $b\in x$ and $y=x\setminus \{b\}$. Then $y<x$ if and only if $b=i_{u}$ with $u\in [2k-n]$.
\end{prop}

\textbf{Proof.} Let $C$ be the $\sim$ equivalence class of $x$. The elements of $C$ not larger than $x$ have the form $1(x)\cup\{i_{u+1},...,i_{t}\}$ with $u=0,...,t$. As we have ${|1(x)|=k-t}$, the pair of $x$ is $p(x)=1(x)\cup \{i_{2k-n+1},...,i_{t}\}$.

Suppose that $y<x$. As $p(x)\subset y$, we must have $b=i_{u}$ with some $u\in [2k-n]$.

Now we only need to show that in this case, $p(x)\subset p(y)$ holds as well. If ${b=i_{1}}$, then $p(y),p(x),x,y$ are all elements of $C$, hence we are done.

 Now suppose that $b=i_{u}$ with $u\in \{2,...,2k-n\}$. Look at the interval ${J=[i_{u-1}+1,i_{u}-1]}$. This interval is disjoint from $*(x)$, so it is the union of disjoint intervals $[a_{1},b_{1}],...,[a_{r},b_{r}]\in I(x)$. Hence,
 $$c_{x}(i_{u-1},i_{u})=2+\sum_{j=1}^{r} c_{x}(a_{j},b_{j})=2.$$ But then $c_{y}(i_{u-1},i_{u})=0$, so $sg(y)_{i_{u-1}}=1$ and $sg(y)_{i_{u}}=0$. Also, $sg(x)$ and $sg(y)$ agrees on every coordinate other than $i_{u-1}$ and $i_{u}$. Thus, $$y=1(y)\cup \{i_{1},...,i_{t}\}\setminus \{i_{u}\}$$ and $$p(y)=1(y)\cup \{i_{2k-n+1},...,i_{t}\}=p(x)\cup \{i_{u-1}\}.$$ This shows that $p(x)\subset p(y)$. \hfill$\Box$

\begin{figure}
 \includegraphics[scale=2]{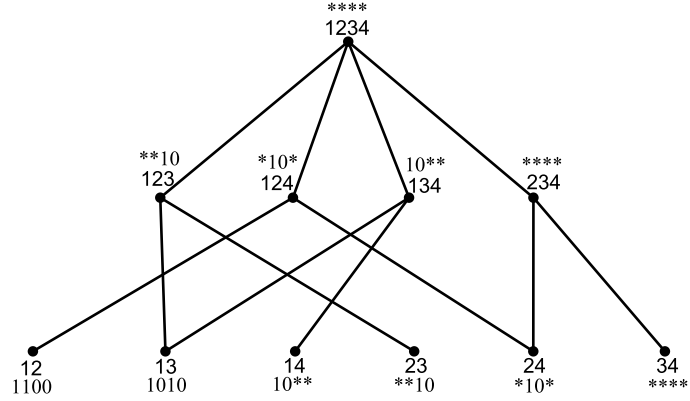}
 \caption{The poset $(Q,<)$ for $n=4$.}
 \label{image1}
 \end{figure}

\bigskip

By Proposition \ref{comparelemma}, the bipartite comparability graph of $(Q,<)$ induced on the levels $[n]^{(k-1)}$ and $[n]^{(k)}$ is $2k-n$ regular from $[n]^{(k)}$. However, it is not regular from $[n]^{(k-1)}$. Also, while this proposition lets us identify the neighbours of any $x\in [n]^{(k)}$ in $[n]^{(k-1)}$ easily, finding a description of the neighbours of $y\in [n]^{(k-1)}$ in $[n]^{(k)}$ is more troublesome.

\bigskip

However, we shall overcome this obstacle by slightly modifying the definition of the signature. We introduce the \emph{circular signature} of an element $x\in 2^{[n]}$. The circular signature of $x$ is denoted by $csg(x)$, and is an element of $\{0,1,*\}^{n}$ defined as follows.

View $x$ as a subset of $\mathbb{Z}_{n}$, the ring of integers modulo $n$. For $i,j\in \mathbb{Z}_{n}$, define the interval $[i,j]_{n}\subset\mathbb{Z}_{n}$ such that
$$[i,j]_{n} = \begin{cases} [i,j] &\mbox{if } i\leq j \\
[i,n]\cup [1,j] & \mbox{if } i>j. \end{cases}$$

For $x\subset 2^{[n]}$, define the function $$c'_{x}(i,j)=|[i,j]_{n}\cap x|-|[i,j]_{n}\setminus x|.$$ If $i\in x$ and there exists $j\in \mathbb{Z}_{n}$ such that $c'_{x}(i,j)=0$, then $csg(x)_{i}=1$ and the \emph{circular pair of $i$ in $x$} is the $j$ with this property for which $|[i,j]_{n}|$ is the smallest and is denoted by $pr'_{x}(i)$. If there is no such $j$, then $csg(x)_{i}=*$.

If $i\not\in x$ and there exists $j\in \mathbb{Z}_{n}$ such that $c_{x}(j,i)=0$, then $csg(x)_{i}=0$ and the \emph{circular pair of $i$ in $x$} is the $j$ with this property for which $|[j,i]_{n}|$ is the smallest and is denoted by $pr'_{x}(i)$. If there is no such $j$, then $csg(x)_{i}=*$.

Also, let $*'(x)=\{i\in [n]: csg(x)_{i}=*\}$ and define $0'(x)$ and $1'(x)$ similarly. Finally, let $I'(x)=\{[i,pr_{x}(i)]_{n}: i\in 1'(x)\}$.

\bigskip

We note that the idea of circular signature can be also found in \cite{necklace}. Taking our earlier example $x=\{2,3,4,6,7,10,11\}\in 2^{[12]}$, now we have $$csg(x)=(\overline{0},*,*,\overline{1,0},\overline{1,\overline{1,0},0},\overline{1,\overline{1,0}}).$$ Also, $*'(x)=\{2,3\}$, $0'(x)=\{1,5,8,9,12\}$ and $1'(x)=\{4,6,7,10,11\}$.

The following proposition lists the properties of the circular signature and helps us compare it with the signature.

\begin{prop}\label{properties_of_csg}
Let $x\in [n]^{(k)}$ with $ n/2\leq k\leq n$. Let $i_{1}<...<i_{t}$ be the elements of $*(x)\cap x$ and $j_{1}<...<j_{t-2k+n}$ be the elements of $*(x)\setminus x$.
\begin{description}
   \item[(i)] If $i\in [n]$ and $j=pr'_{x}(i)$ exists, then $pr'_{x}(j)$ exists and $i=pr'_{x}(j)$. Also, $csg(x)_{i}\neq csg(x)_{j}$.
   \item[(ii)] If $i\in [n]\setminus *(x)$, then $csg(x)_{i}=sg(x)_{i}$ and $pr_{x}(i)=pr'_{x}(i)$.
   \item[(iii)] Any two intervals in $I'(x)$ are either disjoint or one is contained in the other. Also, the intervals are disjoint from $*'(x)$.
   \item[(iv)] If $1\leq r\leq 2k-n$, then $csg(x)_{i_{r}}=*$.
   \item[(v)] If $2k-n<r\leq t$, then $csg(x)_{i_{r}}=1$ and $pr'_{x}(i_{r})=j_{t+1-r}$. Also, if $1\leq s\leq t-2k+n$, then $csg(x)_{j_{s}}=0$ and $pr'_{x}(j_{s})=i_{t+1-s}$.
\end{description}
\end{prop}

\textbf{Proof.} (i) Look at case $i\in x$, the proof in the other case is similar. The function $c'_{x}(i,y)$ changes by 1 as $y$ changes by $1 \mod n$, and $c'_{x}(i,i)=1$. Hence, $c'_{x}(i,y)>0$ for all $y\in [i,j-1]_{n}$, and so we must have $j\not\in x$, proving that $csg(x)_{i}=1-csg(x)_{j}$.

Also, we have $c'_{x}(i,j)=0$, so $pr'_{x}(j)$ exists. Let $i'\in [n]$ be such that $c'_{x}(i',j)=0$ and $|[i',j]_{n}|$ is minimal. Then $i'\in [i,j]$. If $i\neq i'$, then $c'_{x}(i,i')=c'_{x}(i,j)-c'_{x}(i',j)=0$ and $|[i,i']_{n}|<|[i,j]_{n}|$, which is a contradiction.

(ii) If $i<j$, then $c_{x}(i,j)=c'_{x}(i,j)$, so this part trivially follows.

(iii) Let $[a,b]_{n}\in I'(x)$.  The function $c'_{x}(a,y)$ changes by $1$ as $y$ changes by $1\mod n$. Also, $c'_{x}(a,a)=1$, so $c'_{x}(a,y)>0$ for all $y\in [a,b]_{n}$, $y\neq b$. Let $c\in [a,b]_{n}$, $c\neq a,b$. Then $c'_{x}(c,b)<0$. But $c'_{x}(c,c)=1$, hence there exist $d\in [c,b-1]_{n}$ such that $c'_{x}(c,d)=0$. This means that $[a,b]_{n}$ is disjoint from $*'(x)$, and every interval in $I'(x)$ intersecting $[a,b]_{n}$ is either contained in $[a,b]_{n}$ or contains $[a,b]_{n}$.

(iv)-(v) If $c'_{x}(i_{r},j)=0$ for some $j\in [n]$, then the $j$ for which $|[i_{r},j]_{n}|$ is minimal, must be equal to one of $j_{1},...,j_{t-2k+n}$. Suppose $j=j_{r'}$. As $j_{r'}<i_{r}$, we have $c'_{x}(i_{r},j_{r'})=c_{x}(i_{r},n)+c_{x}(1,j_{r'})$. Let $[a_{1},b_{1}],...,[a_{s},b_{s}]$ be the maximal intervals in $I(x)$ contained in $[i_{r},n]$. Then,
$$c_{x}(i_{r},n)=t-r+1+\sum_{l=1}^{s}c_{x}(a_{l},b_{l})=t-r+1.$$
We can show similarly that $c_{x}(1,j_{r'})=-r'$. Hence, $c'_{x}(i_{r},j_{r'})=0$ if and only if $r'=t-r+1$. \hfill$\Box$

\bigskip

Now we can describe the neighbours of any element with the help of the circular signature.

\begin{prop}\label{propimportant}
\begin{description}
\item{(i)} Let $x\in [n]^{(k)}$ with $ n/2\leq k\leq n$ and let $i_{1}<...<i_{2k-n}$ be the elements of $*'(x)$. The elements in $[n]^{(k-1)}$ that are $<$-smaller than $x$  are $x\setminus \{i_{1}\},...,x\setminus \{i_{2k-n}\}$. Furthermore, for $u=1,...,2k-n$, the interval $[i_{u-1},i_{u}]_{n}$ is a maximal interval in $I'(x\setminus \{i_{u}\})$, where the index $u-1$ meant modulo $2k-n$.

\item{(ii)} Also, if $y\in [n]^{(k-1)}$ and $[j_{1},pr_{y}(j_{1})]_{n},...,[j_{s},pr_{y}(j_{s})]_{n}$ are the maximal intervals in $I'(y)$, then there are $s$ elements in $[n]^{(k)}$ that are $<$-larger than $y$, namely $y\cup pr'_{y}(j_{r})$ for $r=1,...,s$.
  \end{description}
\end{prop}

\textbf{Proof.} (i) The first part of (i) is straightforward from Proposition \ref{comparelemma} and Proposition \ref{properties_of_csg}. Now let $z=x\setminus \{i_{u}\}$. We have
$$c'_{z}(i,j) = \begin{cases} c'_{x}(i,j) &\mbox{if } i_{u}\not\in[i,j]_{n} \\
c'_{x}(i,j)-2 & \mbox{if } i_{u}\in[i,j]_{n}. \end{cases}$$

If $i\in 1'(x)$, then $[i,pr'_{x}(i)]_{n}$ does not contain $i_{u}$, as every interval in $I'(x)$ is disjoint from $*'(x)$. Hence, $i\in 1'(z)$ as well, and $pr'_{z}(i)=pr'_{x}(i)$. Also, let $[a_{1},b_{1}]_{n},...,[a_{s},b_{s}]_{n}$ be the maximal intervals in $I'(x)$ contained in $[i_{u-1},i_{u}]_{n}$. Then,   $$c'_{z}(i_{u-1},i_{u})=\sum_{l=1}^{s}c'_{z}(a_{l},b_{l})=\sum_{l=1}^{s}c'_{x}(a_{l},b_{l})=0.$$
Hence, $i_{u-1}\in 1'(z)$, $i_{u}=0'(z)$ and $pr'_{z}(i_{u-1})=i_{u}$ as every element in $[i_{u-1}+1,i_{u}-1]_{n}$ has a pair in $x$. As $|*'(z)|=2k-n-2$, we must have $*'(z)=*'(x)\setminus\{i_{u-1},i_{u}\}$ and so $I'(z)=I'(x)\cup\{[i_{u-1},i_{u}]_{n}$. The interval $[i_{u-1},i_{u}]_{n}$ is also maximal, as no other interval of $I'(z)$ contains $i_{u}$.

(ii) By the second part of (i), if $y\cup \{k\}>y$ for some $k\in [n]$, then $k\in \{pr'_{y}(j_{1}),...,pr'_{y}(j_{r})\}$. Let $k_{r}=pr'_{y}(j_{r})$ and $w=y\cup\{k_{r}\}$. Then, we have $$c'_{w}(i,j) = \begin{cases} c'_{y}(i,j) &\mbox{if } k_{r}\not\in[i,j]_{n} \\
c'_{y}(i,j)+2 & \mbox{if } k_{r}\in[i,j]_{n}. \end{cases}$$
Let $a\in 1'(y)\setminus \{j_{r}\}$. As $[j_{r},k_{r}]$ is a maximal interval in $I'(y)$, we have ${k_{r}\not\in [a,pr'_{y}(a)]_{n}}$. Hence, $c'_{w}(a,b)=c'_{y}(a,b)$ for all $b\in [a,pr'_{y}(a)]_{n}$.
But this means that $a\in 1'(w)$ and $pr'_{w}(a)=pr'_{y}(a)$. Similarly, if $b\in 0'(y)\setminus\{k_{r}\}$, then $b\in 0'(w)$ and $pr'_{w}(b)=pr'_{y}(b)$.

 The sets $1'(w)$ and $0'(w)$ have exactly $n-k$ elements, so we must have $1'(w)=1'(y)\setminus \{j_{r}\}$, $0'(w)=0'(y)\setminus \{k_{r}\}$ and $*'(w)=*'(y)\cup \{j_{r},k_{r}\}$. Hence, by the first part of (i), we have $y<w$. \hfill$\Box$

\bigskip

We finished analyzing the poset $(Q,<)$. But before we can start the proof of Theorem \ref{Qnormalized}, we still need the following well known properties of normalized matching graphs and posets.

\begin{prop}\label{levelsprop}
Let $(P,<)$ be a graded poset with levels $A_{0},...,A_{n}$. For $k=0,...,n-1$, let $G_{k}=(A_{k},A_{k+1},E_{k})$ be the bipartite graph, where $x\in A_{k}$ and $y\in A_{k+1}$ are joined by an edge if $x<y$. If $G_{k}$ is a normalized matching graph for $k=0,...,n-1$, then $(P,<)$ is a normalized matching poset.
\end{prop}

\textbf{Proof.} We need to show that for any positive integers $i$ and $j$ with $0\leq i,j\leq n$, $i\neq j$, the bipartite subgraph induced on $A_{i}\cup A_{j}$ is a normalized matching graph. Suppose that $i<j$, the other case being similar. Let $X_{0}\subset A_{i}$ and for $l=1,...,j-i$, define $X_{l}$ to be the set of elements of $A_{i+l}$ which are larger than some element of $X_{l-1}$. We need to show that $$\frac{|X_{0}|}{|A_{i}|}\leq\frac{|X_{j-i}|}{|A_{l}|}.$$ But this is obvious as we have
$$\frac{|X_{l}|}{|A_{i+l}|}\leq \frac{|X_{l+1}|}{|A_{i+l+1}|}$$
for $l=0,...,j-i-1$, because $G_{i+l}$ is a normalized matching.\hfill$\Box$

\bigskip

The following well known result can be found in various sources \cite{daykin, kleitman2}, we state it without proof.

\begin{prop}\label{normmatchprop}
Let $G=(A,B,E)$ be a bipartite graph. The graph $G$ is a normalized matching graph if and only if there exist positive reals $a$ and $b$, and a weight function ${w: E\rightarrow \mathbb{R}^{+}}$ such that for any $v\in A$ we have $$\sum_{\substack{v\in e\\ e\in E}} w(e)=a,$$ and for any $w\in B$ we have $$\sum_{\substack{w\in e\\ e\in E}} w(e)=b.$$
\end{prop}
\hfill$\Box$

\bigskip

Now we are ready to prove Theorem \ref{Qnormalized}, namely that $(Q,<)$ is a normalized matching poset.

\bigskip

\textbf{Proof of Theorem \ref{Qnormalized}.} By Proposition \ref{levelsprop}, it is enough to show that the bipartite subgraph of the comparability graph of $(Q,<)$ induced on ${[n]^{(k-1)}\cup [n]^{(k)}}$ is a normalized matching graph for $k=\lceil n/2\rceil+1,...,n$.

Fix $k$ with $n/2+1\leq k\leq n$, and let $G=([n]^{(k-1)},[n]^{(k)},E)$ be the bipartite graph, where $x\in [n]^{(k-1)}$ and $y\in[n]^{(k)}$ are joined by an edge if $x<y$. We show that with the choice $a=2n-2k+2$ and $b=2k$ there is a weight function $w: E\rightarrow \mathbb{R}^{+}$ satisfying the conditions of Proposition \ref{normmatchprop}.

Define $w: E\rightarrow \mathbb{R}^{+}$ as follows. Suppose that $x\in [n]^{(k-1)}$ and $y\in [n]^{(k)}$ such that $x<y$ and let the elements of $*'(y)$ be $i_{1}<...<i_{2k-n}$. Then $x=y\setminus \{i_{u}\}$ for some $1\leq u\leq 2k-n$. Also, $[i_{u-1},i_{u}]_{n}$ is a maximal interval in $I'(x)$, where the index $u-1$ meant modulo $2k-n$. Let $$w(\{x,y\})=|[i_{u-1},i_{u}]_{n}|.$$

We show that $w$ suffices. Fix $x\in [n]^{(k-1)}$ and let $[j_{1},pr'_{x}(j_{1})]_{n},...,[j_{t},pr'_{x}(j_{t})]_{n}$ be the maximal intervals in $I'(x)$. The neighbours of $x$ in $[n]^{(k)}$ are $$x\cup \{pr'_{x}(j_{1})\},...,x\cup \{pr'_{x}(j_{t})\}.$$
Hence, we have
$$\sum_{\substack{x\in e\\e\in E}}w(e)=\sum_{u=1}^{t}|[j_{u},pr'_{x}(j_{u})]_{n}|.$$
But the intervals $[j_{1},pr'_{x}(j_{1})]_{n},...,[j_{t},pr'_{x}(j_{t})]_{n}$ are pairwise disjoint and their union is $[n]\setminus *'(x)$, so
$$\sum_{u=1}^{t}|[j_{u},pr'_{x}(j_{u})]_{n}|=|[n]\setminus *'(x)|=2n-2k+2=a.$$

Now fix $y\in [n]^{(k)}$ and let $i_{1}<...<i_{2k-n}$ be the elements of $*'(y)$. The neighbours of $y$ in $G$ are
$$y\setminus\{i_{1}\},...,y\setminus\{i_{2k-n}\}.$$
Hence,
$$\sum_{\substack{y\in e\\e\in E}}w(e)=\sum_{u=1}^{2k-n}|[i_{u-1},i_{u}]_{n}|,$$
where the index $u-1$ is taken modulo $2k-n$. Every element of $[n]\setminus*'(x)$ is contained in exactly one of the intervals $[i_{1},i_{2}]_{n},...,[i_{2k-n},i_{1}]_{n}$, while every element of $*'(x)$ is contained in exactly two of those intervals. Hence,
$$\sum_{u=1}^{2k-n}|[i_{u-1},i_{u}]_{n}|=2k=b.$$
So the weight function $w$ suffices and $(Q,<)$ is a normalized matching poset. \hfill$\Box$ $\Box$

\section{Partitioning the Boolean lattice into rank-symmetric chains of uniform size}

In this section, we prove Theorem \ref{mainthm}. In the proof, we apply the following two theorems from \cite{me}.

\begin{theorem}\label{upper}
Let $P$ be an unimodal normalized matching poset of width $w$. Then the poset $P$ can be partitioned into $w$ chains of size at most $\frac{2|P|}{w}+5$.
\end{theorem}

\begin{theorem}\label{lower}
Let $P$ be a normalized matching poset with levels $A_{0},A_{1},\ldots,A_{n}$ and let $a_{i}=|A_{i}|$ for $i=0,\ldots,n$. Suppose that $w=a_{0}\geq a_{1}\geq\ldots\geq a_{n}$. Let  $f:\{0,\ldots,n\}\rightarrow \mathbb{N}\cup\{\infty\}$ be defined by
$$f(k)=\min\{i>k: a_{k+1}+\ldots+a_{i}\geq w\},$$
where $a_{i}=0$ if $i>n$. Also, let $f_{1}=f$ and $f_{i}=f\circ f_{i-1}$ denote the iterations of $f$; set $f_{0}\equiv 0$. Finally, let $d$ be the largest integer such that $f_{d}(0)<\infty$. The poset $P$ can be partitioned into $w$ chains of size at least $d+1$.
\end{theorem}

The proof of Theorem \ref{mainthm} is almost the same as the proof of Theorem 1.5 in \cite{me}. For completeness, we provide the proof here as well, but we shall copy most of it word by word.

\bigskip

\textbf{Proof of Theorem \ref{mainthm}.}  Again, let $Q=[n]^{(\geq n/2)}$ and choose a bijection $p: Q\rightarrow [n]^{(\leq n/2)}$ and partial ordering $<$ such that they satisfy the conditions of Theorem \ref{mainthm2}. The width of $(Q,<)$ is $w=\binom{n}{\lfloor n/2\rfloor}$ as the symmetric chain decomposition of $2^{[n]}$ also defines a chain decomposition of $(Q,<)$ into $w$ chains. (It also follows from the normalized matching property.) For $i=0,...,\lfloor n/2\rfloor$, let $A_{i}=[n]^{(\lceil n/2\rceil+i)}$. First, we need some approximation of the size of the level $A_{i}$, where $i=O(\sqrt{n})$.

 Using Stirling's approximation one can easily get the estimation that for any $t\in \mathbb{R}$ we have
\begin{equation}\label{5th}
\binom{n}{n/2+t\sqrt{n}}=(1+o(1))2^{n}e^{-2t^{2}}\sqrt{\frac{2}{\pi n}}=(1+o(1))e^{-2t^{2}}w,
\end{equation}
as $n\rightarrow \infty$.

For $k=1,2,\ldots$, let $T_{k}>0$ be the smallest integer such that ${\binom{n}{\lceil n/2\rceil+T_{k}}<w/k}$. Then, by (\ref{5th}) we have $T_{k}=(1/\sqrt{2}+o(1))\sqrt{n}\sqrt{\log k}$.

Let $K$ be the smallest positive integer such that
$$\alpha-\frac{\epsilon}{2}<\sqrt{2}\sum_{k=2}^{K}\frac{\sqrt{\log k}-\sqrt{\log (k-1)}}{k}.$$
Note that for $k>2$ we have
$$\frac{\sqrt{\log k}-\sqrt{\log (k-1)}}{k}=\frac{(\log k)-(\log (k-1))}{k(\sqrt{\log (k-1)}+\sqrt{\log k})}<\frac{1}{k^{2}}.$$
Hence, the tail of the sum $\displaystyle\sum_{k=2}^{\infty}\frac{\sqrt{\log k}-\sqrt{\log (k-1)}}{k}$ can be easily bounded:
$$\sum_{k=K+1}^{\infty}\frac{\sqrt{\log k}-\sqrt{\log (k-1)}}{k}<\sum_{k=K+1}^{\infty}\frac{1}{k^{2}}<\frac{1}{K}.$$
This implies $K=O(1/\epsilon)$.

\bigskip

Let $P$ be the subposet of $(Q,<)$ induced by the levels $A_{0},\ldots,A_{T_{K}}$. Define $f$ and $d$ as in Theorem \ref{lower}. For $k=2,3,\ldots, K$, if $T_{k-1}<j\leq T_{k}-k$, then $f(j)=j+k$. Hence,
$$d\geq \sum_{k=2}^{K}\left\lfloor \frac{T_{k}-T_{k-1}-k}{k}\right\rfloor=\left(\frac{1}{\sqrt{2}}+o(1)\right)\sqrt{n}\sum_{k=2}^{K}\frac{\sqrt{\log k}-\sqrt{\log (k-1)}}{k}\geq$$
$$\geq \left(\frac{1}{2}+o(1)\right)\sqrt{n}\left(\alpha-\frac{\epsilon}{2}\right).$$
Thus, by Theorem \ref{lower} there exists a chain partition of $P$ into $w$ chains, each  of size at least $(1/2+o(1))\sqrt{n}(\alpha-\epsilon/2)$ and at most $$T_{K}=\left(\frac{1}{\sqrt{2}}+o(1)\right)\sqrt{n}\sqrt{\log K}.$$ Let $\{C_{x}\}_{x\in B_{\lceil n/2\rceil}}$ be such a partition with $x\in C_{x}$ for all $x\in B_{\lceil n/2\rceil}$.

Let $P'$ be the subposet of $2^{[n]}$ induced by the levels $A_{T_{K}},\ldots,A_{\lfloor n/2\rfloor}$. Then $Q$ is also a unimodal normalized matching poset with width $$w'=|A_{T_{K}}|=\left(\frac{1}{K}+o(1)\right)w.$$ Thus, by Theorem \ref{upper}, $P'$ has a partition into $w'$ chains, each of size at most
$$\frac{2|P'|}{w'}+5=\frac{2(K+o(1))|Q|}{w}<\frac{(K+o(1))2^{n}}{w}=(K+o(1))\sqrt{\frac{\pi n}{2}}.$$
Let $\{D_{y}\}_{y\in B_{\lceil n/2\rceil+T_{K}}}$ be such a partition with $y\in D_{y}$ for all $y\in B_{\lceil n/2\rceil+T_{K}}$.

For all $x\in B_{\lceil n/2\rceil}$, let $C'_{x}=C_{x}$ if $C_{x}\cap A_{T_{K}}=\emptyset$ and let $C'_{x}=C_{x}\cup D_{y}$ if $C_{x}\cap A_{T_{K}}=\{y\}$. Then $\{C'_{x}: x\in A_{0}\}$ is a partition of $Q$ satisfying
$$\left(\frac{1}{2}+o(1)\right)\left(\alpha-\frac{\epsilon}{2}\right)\sqrt{n}< |C_{x}|<T_{K}+(K+o(1))\sqrt{\frac{\pi n}{2}}=$$
$$=\left(\sqrt{\frac{\log K}{2}}+K\sqrt{\frac{\pi}{2}}+o(1)\right)\sqrt{n}=O\left(\frac{\sqrt{n}}{\epsilon}\right)$$
for all $x\in A_{0}$.

But then $\{C'_{x}\cup p(C'_{x}): x\in A_{0}\}$ is a chain partition of $2^{[n]}$ into rank-symmetric chains of size between $(\alpha-\epsilon/2+o(1))\sqrt{n}$ and $O(\sqrt{n}/\epsilon)$.

Setting $\epsilon=0.04$, one can get the exact bounds $0.8\sqrt{n}$ and $13\sqrt{n}$ for sufficiently large $n$. We shall avoid doing these calculations.
\hfill$\Box$

\section{The number of symmetric chain partitions of the Boolean lattice}

In this section, we present bounds on the number of symmetric chain decompositions of $2^{[n]}$ and we prove Theorem \ref{nbofpartitions}. First, we show a short proof of a slightly worse bound than the one we have in Theorem \ref{nbofpartitions}, without using Theorem \ref{mainthm2}. This simple proof uses the elegant idea of Kleitman \cite{kleitman}, who builds a symmetric chain partition of $2^{[n]}$ by induction.

\begin{prop}
The Boolean lattice $2^{[n]}$ has at least $2^{\Omega(2^{n}/\sqrt{n})}$ symmetric chain decompositions.
\end{prop}

\textbf{Proof.} Let $m=\binom{n-1}{\lfloor (n-1)/2\rfloor}$ and let $C_{1},...,C_{m}$ be a symmetric chain decomposition of $2^{[n-1]}$ such that $|C_{1}|\geq ...\geq |C_{m}|$. Let $1\leq l\leq m$ be the largest integer such that $|C_{l}|>1$, then $$l=\binom{n-1}{\lfloor (n-2)/2\rfloor}=\Theta(2^{n}/\sqrt{n}).$$ For each sequence $(i_{1},...,i_{l})\in \{0,1\}^{l}$, we define the following symmetric chain partition of $2^{[n]}$.

For $j=l+1,...,m$, let $C_{j}^{1}=\{x,x\cup\{n\}\}$, where $x$ is the only element of $C_{j}$. For $j\in [l]$, let the elements of $C_{j}$ be $x_{1}\subset...\subset x_{k}$. If $i_{j}=0$, let
$$C_{j}^{1}=\{x_{1},...,x_{k},x_{k}\cup\{n\}\},$$
$$C_{j}^{2}=\{x_{1}\cup\{n\},...,x_{k-1}\cup\{n\}\};$$
 and if $i_{j}=1$, let
$$C_{j}^{1}=\{x_{1},x_{1}\cup \{n\},...,x_{k}\cup\{n\}\},$$
$$C_{j}^{2}=\{x_{2},...,x_{k}\}.$$
Then one can easily check that $C_{1}^{1},...,C_{m}^{1},C_{1}^{2},...C_{l}^{2}$ is a symmetric chain decomposition of $2^{[n]}$, and for each sequence $(i_{1},...,i_{l})\in \{0,1\}^{l}$ we get a different partition. Hence, we showed that $2^{[n]}$ has at least $2^{l}=2^{\Omega(2^{n}/\sqrt{n})}$ symmetric chain decompositions.\hfill $\Box$

\bigskip

Now we show that we can gain az extra $\log n$ factor in the exponent by utilizing Theorem \ref{mainthm2}.

Let $Q=[n]^{(\geq n/2)}$, and choose a bijection $p$ and partial ordering $<$ satisfying the conditions in Theorem \ref{mainthm2}. Also, for simplicity, write $M=\lceil n/2\rceil$ and $w=\binom{n}{M}$. If $C\subset Q$ is a skipless chain with minimal element lying in $[n]^{M}$, then $C\cup p(C)$ is a symmetric chain in $2^{[n]}$.

For $i=M,...,n-1$, let $G_{i}$ be the comparability graph of the subposet of $(Q,<)$ induced on the levels $[n]^{(i)}$ and $[n]^{(i+1)}$. Then $G_{i}$ is a normalized matching graph. Let $\mathfrak{T}_{i}$ be the set of complete matchings from $[n]^{(i+1)}$ to $[n]^{(i)}$ in $G_{i}$. Note that each sequence of matchings $(T_{M},...,T_{n-1})\in \mathfrak{T}_{M}\times...\times \mathfrak{T}_{n-1}$ corresponds to a unique partition of $Q$ into skipless chains with minimal elements in $[n]^{(M)}$. Look at the graph with vertex set $Q$, where two elements are joined by an edge if they are joined by an edge in one of $T_{M},...,T_{n-1}$. Then the components of this graph are skipless chains, and every component intersects $[n]^{(M)}$. Hence, the components of this graph form a chain partition of $Q$ into skipless chains with minimal elements in $[n]^{(M)}$.

Thus, every element of $\mathfrak{T}_{M}\times...\times \mathfrak{T}_{n-1}$ corresponds to a unique symmetric chain partition of $2^{[n]}$. Our task is reduced to estimating the size of $\mathfrak{T}_{i}$ for $i=M,...,n-1$.

\begin{theorem}\label{matchingsnormalized}
Let $G=(A,B,E)$ be a normalized matching graph with $|A|\leq|B|$. The number of complete matchings from $A$ to $B$ is at least $\sqrt{\binom{|B|}{|A|}}$.
\end{theorem}

\textbf{Proof.} Let $|A|=k$ and $|B|=n$. Let $$M(G)=\{V\in B^{(k)}: \exists\ complete\ matching\ from\ A\ to\ V\ in\ G\}.$$
It is enough to show that $M(G)\geq \sqrt{\binom{n}{k}}$.

We use the following observation: let $G'=(A,B,E')$ be another normalized matching graph on the same vertex set. Then $M(G)\cap M(G')\neq \emptyset$. Create the following graded poset $(P,<)$: let the levels of $P$ be $A,B,A'$, where $A'$ is a disjoint copy of $A$; let the partial order $<$ be induced by the following relations: for $a\in A$, $a'\in A'$ and $b\in B$ let $a<b$ if $ab\in E$, and  $b<a'$ if $a'b\in E'$. Then $(P,<)$ is a unimodal, rank-symmetric normalized matching poset, so by the theorem of Griggs \cite{symchain} it has a symmetric chain decomposition. This chain decomposition contains chains of size 1 and 3. Let $V\subset B$ be the set of elements which are the middle elements of chains of size 3. Then $V\in M(G)\cap M(G')$.

 Let $S_{B}$ denote the set of all permutations of $B$. For $\pi\in S_{B}$ define the bipartite graph $G_{\pi}=(A,B,E_{\pi})$ such that $a\in A$ and $b\in B$ are joined by an edge if $a\pi(b)\in E$.

 Clearly, $G_{\pi}$ is a normalized matching graph as $G_{\pi}$ is isomorphic to $G$. Hence, we have that $M(G)\cap M(G_{\pi})\neq \emptyset$. Note that
$$M(G_{\pi})=\{\pi^{-1}(V): V\in M(G)\},$$ so for each $\pi\in S_{B}$ there exists a pair $(V,W)\in M(G)\times M(G)$ such that $V=\pi^{-1}(W)$. But for each pair $(V,W)$ there are exactly $k!(n-k)!$ permutations $\pi\in S_{B}$ such that $V=\pi^{-1}(W)$. So, we must have
$$k!(n-k)!|M(G)|^{2}\geq n!.$$
\hfill $\Box$

\bigskip

\emph{Remark.} The result of Theorem \ref{matchingsnormalized} is not sharp for any $1\leq k\leq n-1$. However, we cannot get any lower bound better than $\binom{n}{k}$. This is true as every minimal normalized matching graph is a forest (see \cite{daykin}), and in this case the number of complete matchings from $A$ to $B$ is equal to $|M(G)|\leq \binom{n}{k}$.

\bigskip

\textbf{Proof of Theorem \ref{nbofpartitions}.} Write $a_{k}$ for $\binom{n}{k}$. We have already established that $2^{[n]}$ has at least $|\mathfrak{T}_{M}|...|\mathfrak{T}_{n-1}|$ symmetric chain decompositions. Also, by Theorem \ref{matchingsnormalized}, we have $|\mathfrak{T}_{k}|\geq \sqrt{\binom{a_{k}}{a_{k+1}}}$ for $k=M,...,n-1$. We shall estimate $\binom{a_{k}}{a_{k-1}}$ for $n/2+\sqrt{n}<k<n/2+2\sqrt{n}$. In this case, we have $a_{k}=\Theta(2^{n}/\sqrt{n})$ and $a_{k}-a_{k+1}=\Theta(2^{n}/n)$. Hence, using the well known lower bound $\binom{a}{b}\geq(a/b)^{b}$ for binomial coefficients, we get

$$\binom{a_{k}}{a_{k+1}}=\binom{a_{k}}{a_{k}-a_{k+1}}\geq \left(\frac{a_{k}}{a_{k}-a_{k+1}}\right)^{a_{k}-a_{k-1}}=$$
$$=\sqrt{n}^{\Theta(2^{n}/n)}=2^{\Theta(2^{n}\log n/n)}.$$
Thus, we have

$$\prod_{k=M}^{n-1}|\mathfrak{T}_{k}|\geq \prod_{n/2+\sqrt{n}<k<n/2+2\sqrt{n}}\sqrt{\binom{a_{k}}{a_{k+1}}}=2^{\Omega(2^{n}\log n/\sqrt{n})}.$$
\hfill$\Box$

\bigskip

We note that $2^{[n]}$ has at most $2^{O(2^{n}\log n)}$ symmetric chain decompositions, so there is only a $\sqrt{n}$ factor gap in the exponent between the lower and upper bound. We can get the upper bound $n^{2^{n}}$ by the following simple observation: every symmetric chain partition corresponds to a unique sequence $(U_{0},...,U_{n-1})$, where $U_{k}$ is a complete matching from $[n]^{(k)}$ to $[n]^{(k+1)}$ for $k<M$, and $U_{k}$ is a complete matching from $U_{k+1}$ to $U_{k}$ for $k\geq M$. But for each $k$, the number of such matchings is less than $n^{\binom{n}{k}}$, as the comparability subgraph of $2^{[n]}$ induced on the levels $[n]^{(k)}\cup [n]^{(k+1)}$ has maximum degree at most $n$. Hence, the number of symmetric chain decompositions of $2^{[n]}$ is less than $\prod_{k=0}^{n}n^{\binom{n}{k}}=n^{2^{n}}$.

\section{Open problems}

In this section, we propose some open problems.

Let us extend some of our earlier definitions to rank-symmetric posets. Let $(P,<)$ be a rank-symmetric poset with levels $A_{0},...,A_{n}$. A chain $C\subset P$ is \emph{rank-symmetric}, if $|C\cap A_{i}|=|C\cap A_{n-i}|$ for $i=0,...,n$; the chain $C$ is \emph{skipless}, if the set $\{j: |A_{j}\cap C|=1\}$ is an interval. The chain $C$ is \emph{symmetric}, if $C$ is skipless and rank-symmetric.

 Griggs \cite{symchain} showed that every rank-symmetric, unimodal normalized matching poset $(P,<)$ has a chain decomposition into symmetric chains.

After the method we used to prove Theorem \ref{mainthm}, it is natural to ask the following question. Let $P$ be a rank-symmetric, unimodal normalized matching poset of width $w$. Can we partition $P$ into $w$ rank-symmetric chains such that the sizes of the chains are as close to each other as possible?

As we do not know the answer in the simpler case when our chains are not necessarily rank-symmetric, we ask a less ambitious question.

\begin{OQ}\label{question}
Define the function $f:\mathbb{Q}^{+}\rightarrow \mathbb{N}$ as follows: for $r\in \mathbb{Q}^{+}$, $f(r)$ is the maximal positive integer such that any rank-symmetric, unimodal normalized matching poset $P$ of width $|P|/r$ can be partitioned into rank-symmetric chains of size at least $f(r)$. Is it true that $\lim_{r\in\mathbb{Q}} f(r)=\infty$?
\end{OQ}

Another closely related question is motivated by the following theorem of Lonc \cite{lonc}. He proved that if given a positive integer $h$ and $n$ is sufficiently large, then $2^{[n]}$ has a chain partition, where all but at most one of the chains have size $h$. The author of this paper \cite{me2} showed that the smallest such $n$ is $O(h^{2})$ and this bound is best possible up to a constant. We propose the following rank-symmetric version of this problem.

\begin{conjecture}\label{symchainconj}
For every positive integer $h$ there exists a positive integer $N(h)$ such that if $n>N(h)$ and $n$ is odd, then $2^{[n]}$ has a partition into rank-symmetric chains, where all but at most one of the chains have size $2h$.
\end{conjecture}

We note the following connection between the results of this paper and Conjecture \ref{symchainconj}. Suppose that the following conjecture is true.

 \begin{conjecture}
 There exists some real $t>0$ such that taking $k=\lfloor n/2+t\sqrt{n}\rfloor$, the bipartite subgraph of $(Q,<)$ induced on $[n]^{(k)}$ and $[n]^{(k+1)}$ has a spanning tree with maximum degree $o(\sqrt{n})$.
 \end{conjecture}

 Then Conjecture \ref{symchainconj} follows from Theorem \ref{mainthm}. Surprisingly, we have not been able to prove the much weaker statement that the bipartite subgraph of $(Q,<)$ induced on $[n]^{(k)}$ and $[n]^{(k+1)}$ is connected for $k=n/2+\Theta(\sqrt{n})$. We note that if $n$ is even, then the bipartite subgraph of $(Q,<)$ induced on the levels $[n]^{(n/2)}$ and $[n]^{(n/2+1)}$ is not even connected, it is the union of trees of size $n+1$.

\end{document}